\documentclass{amsart}

\usepackage[margin=1in]{geometry}

\usepackage{filecontents}
\usepackage{fancyhdr}
\usepackage{lastpage}
\usepackage{bm}
\usepackage[table,dvipsnames]{xcolor}
\usepackage{pgfplots}
\usepackage{caption}
\usepackage[T1]{fontenc}
\usepackage{tabu}
\usepackage{graphicx}
\usepackage{multirow}
\usepackage{amsmath,amssymb,amsthm}
\usepackage{hyperref}
\usepackage{tikz}
\usepackage{tikz-cd}
\usetikzlibrary{calc,matrix,arrows,decorations.markings}
\usepackage{array}
\usepackage{enumerate}
\usepackage{nicefrac}
\usepackage{listings}
\usepackage{centernot}
\usepackage{mathtools}
\usepackage{stmaryrd}
\usepackage{euscript}
\usepackage{manfnt}
\usepackage{marginnote}
\usepackage{mathrsfs}
\usepackage{calligra}
\usepackage{enumitem}
\usepackage[refpage]{nomencl}
\usepackage{accents}

\DeclareMathAlphabet{\mathpzc}{OT1}{pzc}{m}{it}

\newcommand{\F}{\mathbb{F}}
\newcommand{\Fp}{\F_p}

\newcommand{\Q}{\mathbb{Q}}

\newcommand{\C}{\mathbb{C}}

\newcommand{\Z}{\mathbb{Z}}

\newcommand{\m}{{\mathfrak m}}

\newcommand{\Hom}{\operatorname{Hom}}

\newcommand{\A}{\mathbb{A}}
\renewcommand{\P}{\mathbb{P}}

\newcommand{\chr}{\operatorname{char}}

\newcommand{\Spec}{\operatorname{Spec}}

\renewcommand{\H}{\hat{H}}

\newcommand{\Ext}{\operatorname{Ext}}

\newcommand{\Br}{\operatorname{Br}}
\newcommand{\hBr}{\widehat{\Br}}

\newcommand{\op}{{\rm op}}

\newcommand{\sF}{\mathscr{F}}

\newcommand{\sI}{\mathscr{I}}

\newcommand{\sL}{\mathscr{L}}

\newcommand{\sM}{\mathscr{M}}
\newcommand{\sN}{\mathscr{N}}

\newcommand{\sX}{\mathscr{X}}

\newcommand{\sO}{\mathscr{O}}
\newcommand{\sHom}{\mathscr{H}\text{\kern -.35em{\calligra\Large om}\kern .08em}}

\newcommand{\cU}{\mathcal{U}}

\newcommand{\Pic}{\operatorname{Pic}}
\newcommand{\hPic}{\widehat{\operatorname{Pic}}}

\newcommand{\bG}{\mathbb{G}}

\newcommand{\hbG}{\widehat{\bG}}

\newcommand{\et}{\textrm{\'et}}
\newcommand{\Ann}{\operatorname{Ann}}
\newcommand{\doubles}[1]{\llbracket #1\rrbracket}
\newcommand{\height}{\operatorname{ht}}
\newcommand{\MU}{\mathrm{MU}}
\newcommand{\len}{\operatorname{len}}
\newcommand{\rP}{\mathrm{P}}

\newcommand{\cris}{\mathrm{cr}}

\allowdisplaybreaks

\newtheorem{theorem}{Theorem}[section]

\theoremstyle{definition}
\newtheorem{definition}[theorem]{Definition}

\newtheorem{remark}[theorem]{Remark}
\newtheorem{example}[theorem]{Example}
\theoremstyle{remark}

\numberwithin{equation}{section}

\definecolor{zaffre}{rgb}{0.0, 0.08, 0.66}
\hypersetup{colorlinks=true,allcolors=zaffre}

\renewcommand{\H}{\mathrm{H}}
\newcommand{\BP}{\mathrm{BP}}

\title{Constructing Explicit K3 Spectra}
\author{Oron Y. Propp}
\address{Department of Mathematics, Massachusetts Institute of Technology, Cambridge, MA 02139, USA}
\email{opropp@mit.edu}

\begin{document}

\begin{abstract}
Following an overview of the relevant theory, we construct several explicit examples of height-$3$ K3 spectra.
\end{abstract}

\maketitle

\section{Introduction}

Algebraic topology is broadly concerned with associating algebraic invariants, such as cohomology theories, to topological spaces. The subfield of chromatic homotopy theory, which grew out of work of Quillen in the 1960s, allows one to study a broad class of cohomology theories via certain purely algebraic structures called formal group laws. These often arise in turn from algebro-geometric objects, such as the infinitesimal deformations of an abelian variety. In this note, we explore the relationship between K3 surfaces and cohomology theories, taking as a starting point recent work of Szymik \cite{szymik}. We begin with an introduction to chromatic homotopy theory, following \cite{luriechrom}.

Let $E$ be a multiplicative cohomology theory. We say that $E$ is \emph{complex-orientable} if when $X=\C\rP^\infty$, the Atiyah--Hirzebruch spectral sequence
\begin{equation*}
E_2^{i,j}=H^i(X;E^j(*))\implies E^{i+j}(X)
\end{equation*}
degenerates at its second page. In this case, one obtains a (non-canonical) isomorphism $E^*(\C\rP^\infty)\cong E^*(*)\doubles{t}$ with a power series ring in some generator $t\in E^2(*)$; this choice of $t$ is called a \emph{complex orientation} of $E$. One can then define a Chern class taking values in $E$-cohomology analogously to that for singular cohomology: a line bundle $\sL$ on a space $X$ is classified by a continuous map $f\colon X\to\C\rP^\infty$ such that the pullback $f^*\sO(-1)\cong\sL$, where $\sO(-1)$ denotes the tautological bundle on $\C\rP^\infty$; we may then define $c_1^E(\sL):=f^*t\in E^*(X)$.

Consider the space $\C\rP^\infty\times\C\rP^\infty$; it is the universal example of a space with two complex line bundles, namely $\pi_1^*\sO(-1)$ and $\pi_2^*\sO(-1)$, where $\pi_1$ and $\pi_2$ are the two projections $\C\rP^\infty\times\C\rP^\infty\to\C\rP^\infty$. The Atiyah--Hirzebruch spectral sequence gives an isomorphism $E^*(\C\rP^\infty\times\C\rP^\infty)\cong E^*(*)\doubles{x,y}$, where $x:=\pi_1^*t=c_1^E(\pi_1^*\sO(-1))$ and $y:=\pi_2^*t=c_1^E(\pi_2^*\sO(-1))$. Under this isomorphism, the Chern class $c_1^E(\pi_1^*\sO(-1)\otimes\pi_2^*\sO(-1))$ of the universal \emph{tensor product} of two line bundles is a power series $G_E(x,y)\in E^*(*)\doubles{x,y}$. Thus, by universality, $c_1^E(\sL\otimes\sL')=G_E(c_1^E(\sL),c_1^E(\sL'))$ for any two line bundles $\sL$ and $\sL'$ on any space $X$. For the ordinary Chern class, the power series $f$ is given by $G_E(x,y)=x+y$, which translates into the usual multiplicativity formula. However, in general, this formula need not hold, and $f$ can be quite complicated; nonetheless, $G_E$ must satisfy certain identities, corresponding to basic facts about the tensor product of complex line bundles (note that the trivial line bundle corresponds to $0\in E^*(\C\rP^\infty\times\C\rP^\infty)$):
\begin{equation}
\label{eqn:fgl}
\begin{gathered}
G_E(x,0)=x=G_E(0,x)\\
G_E(x,y)=G_E(y,x)\\
G_E(G_E(x,y),z)=G_E(x,G_E(y,z))
\end{gathered}
\qquad
\begin{gathered}
\text{(Unitality),}\\
\text{(Commutativity),}\\
\text{(Associativity).}
\end{gathered}
\end{equation}
In general, given a commutative ring $R$, we refer to a power series $G(x,y)\in R\doubles{x,y}$ satisfying these three identities as a \emph{formal group law} over $R$. If we ``forget'' the choices of coordinates $x$ and $y$, then $G$ becomes a \emph{formal group} $\bG$ (see \S\ref{sec:3}); thus, we may associate to any complex-orientable cohomology theory $E$ a formal group $\bG_E$ that does not depend on the choice of complex orientation.

To better understand this connection between formal group laws and cohomology theories, we turn to the universal case. Consider the power series $G=\sum_{i,j\ge 0}a_{i,j}x^iy^j\in\Z[\{a_{i,j}\}_{i,j\ge 0}]\doubles{x,y}$, where each $a_{i,j}$ is graded in degree $2(i+j-1)$ (so that, if we regard $x$ and $y$ as being graded in degree $-2$, then $G$ is homogeneous in degree $-2$). Though $G$ is not a formal group law, we can force it to be one by quotienting $\Z[\{a_{i,j}\}_{i,j\ge 0}]$ by each of the relations between the $a_{i,j}$ implied by the identities in the definition of a formal group law. So for instance, unitality implies that $a_{i,0}=a_{0,j}=0$ for all $i,j\ge 0$ and $a_{1,0}=a_{0,1}=1$; commutativity implies that $a_{i,j}=a_{j,i}$ for all $i,j\ge 0$. The power series $G$ over the resultant ring $L$, known as the \emph{Lazard ring}, is then the universal formal group law, in the sense that any formal group law $F$ over a commutative ring $R$ arises from a map $L\to R$ sending each $a_{i,j}$ to the corresponding coefficient of $F$ in $R$. Thus, in the previous paragraph, a complex orientation of $E$ gives rise to a graded ring map $L\to E^*(*)$. At first glance, the structure of $L$ seems extremely complicated; however, Lazard was able to show that $L\cong\Z[t_1,t_2,\ldots]$, where each $t_i$ has degree $2i$ \cite{lazard}.

Recall that this is precisely the coefficient ring $\pi_*\MU$ of the ring spectrum $\MU$ representing complex cobordism (we assume a basic familiarity with spectra; see \cite{adams} for an introduction). Thus, $\MU$ is \emph{even}, i.e., its homotopy groups are trivial in odd degrees, and it is \emph{periodic}, i.e., its multiplication $\MU\wedge\MU\to\MU$ induces isomorphisms
\begin{equation*}
\pi_{2m}\MU\otimes_{\pi_0\MU}\pi_{2n}\MU\xrightarrow{\cong}\pi_{2(m+n)}\MU
\end{equation*}
for all integers $m$ and $n$, so it must be complex-orientable. One might hope that, after choosing a complex-orientation of $\MU$, the induced map $L\to \MU^*(*)$ classifying its formal group law would be an isomorphism of graded rings. This is indeed (canonically) the case, as proven by Quillen \cite{quillen}; thus, $\MU$ is the universal complex-orientable cohomology theory.

Using this observation, we can attempt to reverse the above constructions: given a formal group law $G$ over a commutative ring $R$, define a functor from topological spaces to abelian groups via
\begin{equation}
\label{eqn:cohom}
X\mapsto\MU^*(X)\otimes_{\MU^*(*)}R,
\end{equation}
where the tensor product is taken with respect to the graded ring map $\MU^*(*)\to R$ induced by $G$. This will not in general be a cohomology theory; however, Landweber's exact functor theorem (see \S\ref{sec:7}) gives a purely algebraic criterion on $G$ (called \emph{Landweber exactness}) which, if satisfied, guarantees that it is \cite{landweber}. The Brown representability theorem then shows that (\ref{eqn:cohom}) is representable by an even periodic ring spectrum $E$ with coefficient ring $R$, whose formal group $\bG_E$ becomes $G$ after a suitable choice of coordinate \cite{brown}. In fact, we can easily recover many interesting theories, such as complex $K$-theory, in this manner.

Landweber's result leads to a natural question: from where can we obtain interesting formal group laws satisfying Landweber's criterion? One answer is given by the infinitesimal deformations of an abelian variety, such as an elliptic curve. Cohomology theories (``elliptic spectra'') arising in this manner have already been studied deeply \cite{tmf,lurie}; however, their formal group laws are limited to heights $1$ and $2$ (see \S\ref{sec:3}), corresponding to the cases of ordinary and supersingular elliptic curves. To study higher heights, we turn instead to formal group laws arising from the formal Brauer group of a K3 surface (see \S\ref{sec:4}), which have height between $1$ and $10$ (as well as $\infty$). This has led Szymik to make the following definition, in analogy with the notion of an elliptic spectrum:
\begin{definition}[{\cite[\S 4.2]{szymik}}]
\label{def:k3spectrum}
A \emph{K3 spectrum} is a triple $(E,X,\varphi)$ consisting of an even periodic ring spectrum $E$, a K3 surface $X$ over $\pi_*E$, and an isomorphism $\varphi$ of the formal Brauer group of $X$ with the formal group $\bG_E$ of $E$.
\end{definition}
Szymik then shows how to obtain K3 spectra from the moduli stack $\sM_{2d}$ of polarized K3 surfaces of degree $2d$ (i.e., with a choice of ample line bundle whose self-intersection number is $2d$) using the Landweber exact functor theorem:
\begin{theorem}[{\cite[Thm.~1]{szymik}}]
Let $R$ be a Noetherian local $\Z_{(p)}$-algebra for some prime $p$ which does not divide $2d$. Let $X$ be a polarized K3 surface of degree $2d$ over $R$ such that the height of the formal Brauer group of the fiber of $X$ over the closed point of $R$ is finite. If the map $\Spec R\to\sM_{2d}$ classifying $X$ is flat, then the formal Brauer group $\hBr_X$ is Landweber exact, so that there is an even periodic ring spectrum $E$ with $\pi_0E\cong R$ and $\bG_E\cong\hBr_X$.
\end{theorem}
\noindent Szymik goes on to note as consequences that any K3 surface $X$ classified by a geometric point of $\sM_{2d}$ gives rise to a K3 spectrum with coefficient ring isomorphic to the \'etale local ring of $\sM_{2d}$ at $X$, and likewise for the ring of formal functions on the formal deformation space of $X$ in $\sM_{2d}$. Thus, K3 spectra exist in relative abundance.

However, Szymik does not provide any explicit examples of K3 spectra; these are of interest to topologists because they facilitate a more direct examination of cohomology theories of height between $3$ and $10$ (and the geometry of the corresponding locus on the moduli stack of formal groups). Our main goal in this note is to give several such examples (see \S\ref{sec:8}), focusing on height $3$ for simplicity. More broadly, it is intended as a primer on constructing K3 spectra directly via the Landweber exact functor theorem, rather than Szymik's results; throughout, we emphasize explicit examples and computations, while also providing a sense of the general theory.

In \S\ref{sec:2}, we define K3 surfaces and discuss several important examples, as well as some of the theory of elliptic fibrations. In \S\ref{sec:3}, we introduce some basics on formal groups. In \S\ref{sec:4}, we recall some elementary deformation theory, which we use in \S\ref{sec:4n} to define the formal Brauer group of a K3 surface. In \S\ref{sec:5}, we introduce two algorithms of Stienstra for computing the formal Brauer group of a K3 surface described as either a complete intersection in projective space or a double cover of projective space. In \S\ref{sec:6}, we do the same for Artin's algorithm for computing the formal Brauer group of an elliptic K3 surface. In both this section and the last, proofs are included essentially in full; though long and involved, these are in both cases quite elegant. The reader primarily concerned with applications may freely skip these proofs. In \S\ref{sec:7}, we state the Landweber exact functor theorem and work through several examples. Finally, in \S\ref{sec:8}, we construct three examples of height-$3$ K3 spectra, arising from families of K3 surfaces of the different types discussed in \S\ref{sec:2}. We conclude in \S\ref{sec:9} by noting some questions for further study suggested by the constructions in \S\ref{sec:8}.

\subsection*{Acknowledgments} The author wishes to thank Hood Chatham for supervising this research, for his help in learning much of this material, and for his many suggestions and insights; Haynes Miller for proposing this project, and for his valuable ideas and guidance; Atticus Christensen, David Corwin, Christian Liedtke, Davesh Maulik, Jack Morava, Sam Raskin, Catherine Ray, Matthias Sch\"utt, Drew Sutherland, and Jared Weinstein for many fruitful discussions; and Jan Stienstra and Markus Szymik for their comments on an early draft of this paper. This research was generously supported by the \textsc{mit} Mathematics Department's \textsc{urop+} program and the Paul E. Gray (1954) Endowed Fund for \textsc{urop}.

\section{K3 Surfaces}
\label{sec:2}

We begin by introducing our primary algebro-geometric object of study, following \cite{k3} (an excellent reference on this subject):

\begin{definition}
A \emph{K3 surface} over a field $k$ is a smooth proper $k$-variety $X$ of dimension $2$ such that $\Omega_{X/k}^2\cong\sO_X$ and $H^1(X,\sO_X)=0$.
\end{definition}

In the rest of this section, we present examples of K3 surfaces which will later become useful to us in constructing K3 spectra.

\begin{example}
\label{ex:k3}
(i) Any smooth quartic $X\subset\P^3$ is a K3 surface. Indeed, applying the adjunction formula to the conormal exact sequence
\begin{equation*}
0\to\sO_{\P_k^3}(-4)|_X\to\Omega_{\P^3/k}|_X\to\Omega_{X/k}\to 0
\end{equation*}
yields an isomorphism $\Omega_{X/k}^2\cong\Omega_{\P^3/k}|_X^2\otimes\sO_{\P_k^3}(4)|_X\cong\sO_X$. Applying the long exact sequence on cohomology to the closed subscheme exact sequence
\begin{equation*}
0\to\sO_{\P_k^3}(-4)\to\sO_{\P^3}\to\sO_X\to 0
\end{equation*}
then gives the vanishing $H^1(X,\sO_X)=0$, as $H^1(\P^3,\sO_X)=H^2(\P^3,\sO_{\P^3}(-4))=0$. Analogously, any smooth cubic in $\P^2$ is an elliptic curve; however, unlike the case of elliptic curves, not every K3 surface admits an embedding into $\P^3$.

(ii) More generally, if $n\ge 3$, a smooth complete intersection in $\P^n$ cut out by homogeneous polynomials $(f_1,\ldots,f_{n-2})$ is a K3 surface if and only if $\sum_{i=1}^{n-2}\deg f_i=n+1$. The case $n=3$ is simply (1); taking $n=4$ or $n=5$ yields K3 surfaces of degree $6$ or $8$, respectively.

(iii) Consider a double cover $\pi\colon X\to\P^2$ branched along a smooth degree-$6$ plane curve $C$. The Hurwitz formula gives $\Omega_{X/k}^2\cong\pi^*(\Omega^2_{\P^2/k}\otimes\sO_{\P^2}(6)^{\otimes 1/2})\cong\sO_X$, and since $\pi_*\sO_X\cong\sO_{\P^2}\oplus\sO_{\P^2}(-3)$ (see for instance \cite[Ch. 0.\S1]{enriques}), we may deduce $H^1(X,\sO_X)=0$. Thus $X$ is a K3 surface, called a \emph{double plane}, if it is smooth; one readily sees that if $\chr(k)\ne 2$, then it suffices for $C$ to be smooth.
\end{example}

\begin{definition}
An \emph{elliptic K3 surface} is a K3 surface $X$ equipped with a surjective morphism $\pi\colon X\to\P^1$ whose generic geometric fiber is a smooth integral curve of genus one.
\end{definition}

We shall be especially interested in elliptic K3 surfaces $X$ whose elliptic fibration $\pi\colon X\to\P^1$ admits a section, i.e., a curve $C_0\subset X$ such that $\pi|_{C_0}\colon C_0\xrightarrow{\cong}\P^1$ (thus, $C_0$ meets every fiber of $\pi$ transversally, in particular, at a smooth point). One can show in this case that the linear system $\sO(3C_0)|_{X_t}$ on the fibers $X_t$ of $\pi$ defines a morphism
\begin{equation*}
X\to\P(\sO_{\P^1}(4)\oplus\sO_{\P^1}(6)\oplus\sO_{\P^1}),
\end{equation*}
over $\P^1$, whose image $X'$ (the \emph{``Weierstrass model''} of $X$) is described by a Weierstrass equation
\begin{equation*}
y^2z+a_1xyz+a_3yz^2=x^3+a_2x^2z+a_4xz^2+a_6z^3,
\end{equation*}
where $x$, $y$, and $z$ are the local coordinates of the direct summands $\sO_{\P^1}(4)$, $\sO_{\P^1}(6)$, and $\sO_{\P^1}$, respectively, and $a_i\in H^0(\P^1,\sO_{\P^1}(2i))$. The zero-section $C_0$ of $X$ then corresponds to that of $X'$ given by the point $[x:y:z]=[0:1:0]$.

The discriminant of $X'$ is the section $\Delta\in H^0(\P^1,\sO_{\P^1}(24))$ computed in terms of the $a_i$ as for an elliptic curve in Weierstrass form. This shows that $\pi$ has $24$ singular fibers, i.e., over the points where $\Delta$ vanishes, if counted with appropriate multiplicities; the Kodaira classifications of these fibers can be determined via Tate's algorithm \cite{tate}.

Dehomogenizing, we obtain the equation
\begin{equation}
\label{eqn:weierstrass}
y^2+a_1(t)xy+a_3(t)y=x^3+a_2(t)x^2+a_4(t)x+a_6(t),
\end{equation}
where $a_i(t)\in k(t)$, describing $X'$ over an affine patch of $\P^1$. This equation is \emph{minimal} if $v_t(a_i)\ge 0$ for all $i$ and $v_t(\Delta)$ is minimal over all such equations for the Weierstrass model of $X$ (which is always satisfied if $v_t(\Delta)<12$); here $v_t$ denotes the $t$-adic valuation on $k(t)$. In this case, we will have $\deg a_i(t)\le 2i$ for each $i$ and $\deg a_i(t)>i$ for some $i$; conversely, if this condition is satisfied for some minimal Weierstrass equation, then the minimal desingularization of the associated surface is an (elliptic) K3 surface (see for instance \cite[\S II.4]{beukers}). We will regularly use this fact to construct Weierstrass models associated to K3 surfaces.

\begin{example}
\label{ex:fermat}
The \emph{Fermat quartic} $X\subset\P_k^3$ cut out by
\begin{equation*}
x_0^4+x_1^4+x_2^4+x_3^4
\end{equation*}
is smooth as long as $\chr(k)\ne 2$, and is then a K3 surface by Example~\ref{ex:k3}(i). In addition, it admits many elliptic fibrations: after choosing any line $\ell\subset X$, the complete linear system $|\sO(1)\otimes\sO(-\ell)|$ consisting of all hyperplanes in $\P^3$ containing $\ell$ defines a morphism $X\to\P^1$. That is, we may project $X$ away from $\ell$ onto a disjoint line in $\P^3$; the fibers are then the plane cubics obtained by removing $\ell$ from the hyperplane intersections of $X$ containing $\ell$. For instance, take $k=\F_{25}$, and let $\ell$ be the line defined by $x_0=\zeta_8x_1$ and $x_2=\zeta_8x_3$, where $\zeta_8\in\F_{25}$ is a primitive root of unity. Then $(x_1-\zeta_8x_1)+t(x_2-\zeta_8x_3)$ is a family of hyperplanes passing through $\ell$ parametrized by an affine patch of $\P^1$; intersecting such a hyperplane with the Fermat quartic gives the vanishing locus of $(\zeta_8x_1+t(x_2-\zeta_8x_3))^4+x_1^4+x_2^4+x_3^4)$ in $\P^2_{k[t]}$, which has two irreducible components: a line corresponding to $\ell$ and a plane cubic. Computing the minimal model of the latter (e.g., via Magma's \texttt{MinimalModel()} function \cite{magma}) gives a Weierstrass model
\begin{equation}
\label{eqn:wmfq}
y^2 = x^3 + 3t^2x^2 + (4t^{10} + 3t^6 + 4t^2)
\end{equation}
for the Fermat quartic. Note that the coefficients $a_i(t)$ of this equation satisfy the condition in the preceding discussion: we have $\deg a_i(t)\le 2i$ for each $i$ and $\deg a_6(t)=10>6$. This Weierstrass equation has discriminant
\begin{equation*}
\Delta=t^4(t+1)^2(t+2)^2(t+3)^2(t+4)^2(t^2+2)^2(t^2+3)^2,
\end{equation*}
which is indeed minimal as $v_t(\Delta)=4$. Knowing $\Delta$ allows us to deduce the singular fibers of this elliptic fibration: for instance, its fiber over $t=0$ has Kodaira symbol $\mathrm{IV}$, that is, it comprises three curves meeting transversally at a point; its fiber over $t=-1$ has Kodaira symbol $\mathrm{I}_2$, that is, it comprises two curves meeting transversally at two points with multiplicity one.
\end{example}

\section{Formal Groups}
\label{sec:3}

In this section, we introduce some basic notions in the theory of formal groups. Let $R$ be a ring. Denote by $\mathsf{Adic}_R$ the category of \emph{adic $R$-algebras}: its objects are augmented $R$-algebras $\epsilon\colon A\to R$ that are complete with respect to the augmentation ideal $\m_A:=\ker\epsilon$, and its morphisms are $R$-algebra maps commuting with the augmentations. The \emph{formal affine line} is the functor $\widehat{\A}_R^1\colon\mathsf{Adic}_R\to\mathsf{Set}$ sending an adic $R$-algebra $A$ to its augmentation ideal $\m_A$; it is represented by the power-series ring $R\doubles{t}$. A ($1$-dimensional, commutative) \emph{formal group} $\bG$ is a commutative group object in the category of functors $\mathsf{Adic}_R\to\mathsf{Set}$ (equivalently, a functor $\bG\colon\mathsf{Adic}_R\to\mathsf{Ab}$) admitting an isomorphism $\bG\xrightarrow{\cong}\widehat{\A}_R^1$, called a \emph{coordinatization of $\bG$}.

Under such a coordinatization, the group operation $\bG\times\bG\to\bG$ is, by the Yoneda lemma, the same as a map $R\doubles{t}\to R\doubles{x,y}$ of adic $R$-algebras, which is the same as a power series $G(x,y)\in R\doubles{x,y}$. The group axioms for $\bG$ then translate directly into the identities (\ref{eqn:fgl}), which characterize $G$ as a \emph{formal group law} (for $\bG$). We often denote $G(x,y)$ by $x+_Gy$, or leaving the choice of formal group law implicit, $x+_{\bG}y$.

Likewise, if $G_1$ and $G_2$ are formal group laws for formal groups $\bG_1$ and $\bG_2$, respectively, then a morphism $\bG_1\to\bG_2$ of formal groups is the same as a map $R\doubles{t}\to R\doubles{t}$ of adic $R$-algebras, i.e., a power series $f\in R\doubles{t}$, satisfying $f(0)=0$ and
\begin{equation}
\label{eqn:fglhom}
G_2(f(x),f(y))=f(G_1(x,y)).
\end{equation}
Such a power series is called a \emph{homomorphism} $f\colon G_1\to G_2$ of formal group laws. It is an \emph{isomorphism} if there exists a homomorphism $g\colon G_2\to G_1$ such that $f(g(t))=t=g(f(t))$; such a composition inverse exists whenever $f'(0)$ is a unit in $R$.

One homomorphism of particular interest for any formal group $\bG$ is its \emph{$p$-series} $[p]_\bG\colon\bG\to\bG$, defined by
\begin{equation*}
\underbrace{t+_\bG t+_\bG\cdots+_{\bG}t}_{p\text{ times}}\in R\doubles{t}
\end{equation*}
for any prime $p$ (after a choice of coordinate $t$). If $R$ is a field of characteristic $p$, then one can show that the first nonzero term in the power-series expansion of $[p]_\bG$ (in fact, any homomorphism of formal group laws) must either be $0$ or of the form $ut^{p^h}$ for some unit $u\in R^\times$ and integer $h\ge 0$. In the latter case, the \emph{height} $\height(\bG)$ of $\bG$ is defined to be $h$; in the former, we set $\height(\bG):=\infty$. The height is an important invariant of a formal group: one can show that it is independent of the choice of coordinate on $\bG$, hence respects isomorphisms. By a well-known result of Dieudonn\'e and Lubin--Tate, formal groups of every possible height exist over any field; if in addition the field is algebraically closed, then any two formal groups of the same height are isomorphic (see for instance \cite[Thm.~13.8]{coctalos}).

\begin{example}
\label{ex:fg}
(i) The \emph{additive formal group $\widehat{\bG}_a$} is the functor $\widehat{\bG}_a(A):=\m_A$, with the group structure on $\m_A$ given by its usual addition. The identity natural transformation $\widehat{\bG}_a\to\widehat{\A}_R^1$ is a coordinatization of $\widehat{\bG}_a$, under which $x+y\in R\doubles{x,y}$ is a formal group law for $\widehat{\bG}_a$. Its $p$-series is then given by $px$, which is $0$ over any field of characteristic $p$, hence $\height(\widehat{\bG}_a)=\infty$.

(ii) The \emph{multiplicative formal group $\widehat{\bG}_m$} is the functor $\widehat{\bG}_m(A):=(1+\m_A)^\times$, i.e., sending $A$ to the group (under multiplication) of units in $A$ of the form $1+a$ for some $a\in\m_A$. One possible coordinate is the map $(1+\m_A)^\times\to\m_A$ given by $1-a\mapsto a$, under which $\widehat{\bG}_m$ has formal group law $1-(1-x)(1-y)=x+y-xy\in R\doubles{x,y}$. Its $p$-series is then given by $1-(1-t)^p\equiv\pm t^p\bmod p$, hence $\height(\widehat{\bG}_m)=1$.
\end{example}

Now suppose that $R$ is a $\Q$-algebra. We claim that for any formal group $\bG$, there exists a unique isomorphism $\log_{\bG}\colon\bG\to\bG_a$, called the \emph{logarithm} of $\bG$. After choosing a formal group law for $\bG$, such a homomorphism must satisfy
\begin{equation*}
\log_{\bG}(\bG(x,y))=\log_{\bG}(x)+\log_{\bG}(y)
\end{equation*}
by (\ref{eqn:fglhom}). Differentiating with respect to $y$ and setting $y=0$ gives
\begin{equation*}
\log_{\bG}'(x)\frac{\partial\bG}{\partial y}(x,0)=\log_{\bG}'(0),
\end{equation*}
which is a unit in $R$; after a change of coordinate, we may assume that it is $1$. By unitality, $\frac{\partial G}{\partial y}(x,0)$ has constant term $1$, so we may invert it in $R\doubles{x}$ to solve for $\log_{\bG}'(x)$. Since $R$ is a $\Q$-algebra, we may express the logarithm of $\bG$ via the formal integral
\begin{equation}
\label{eqn:log}
\log_{\bG}(x)=\int\!\frac{dx}{\frac{\partial\bG}{\partial y}(x,0)},
\end{equation}
verifying both existence and uniqueness. Inverting the logarithm gives the \emph{exponential} $\exp_{\bG}\colon\bG_a\to\bG$, which allows us to recover the formal group law of $\bG$ from its logarithm via
\begin{equation}
\label{eqn:fglfromlog}
\bG(x,y)=\exp_{\bG}(\log_{\bG}(x)+\log_{\bG}(y)),
\end{equation}
or its $p$-series via
\begin{equation*}
[p]_{\bG}=\exp_{\bG}(p\log_{\bG}(t)).
\end{equation*}

\begin{example}
\label{ex:log}
Under the coordinatizations given in Example~\ref{ex:fg}, both the logarithm and exponential of $\bG_a$ are evidently given by the series $t\in R\doubles{t}$. For $\bG_m$, computing explicitly via (\ref{eqn:log}) gives
\begin{equation*}
\log_{\bG_m}(t)=\int\frac{dt}{1-t}=\int\sum_{m\ge 0}t^mdt=\sum_{m\ge 1}\frac{t^m}{m},
\end{equation*}
and therefore
\begin{equation*}
\exp_{\bG_m}(t)=-\sum_{m\ge 1}(-1)^m\frac{t^m}{m}.
\end{equation*}
\end{example}

For more on formal groups, see \cite{ando,coctalos,haynes,cobordism}.

\section{Deformation Theory}
\label{sec:4}

In this section, we recall some elementary deformation theory, following \cite{schlessinger}, \cite[Ch.~6]{fgaexplained}, and \cite[Ch.~18]{k3}. Let $k$ be a field, and denote by $\mathsf{Art}_k$ the category of local Artin $k$-algebras. A \emph{deformation functor} is a functor $F\colon\mathsf{Art}_k\to\mathsf{Set}$ such that $F(k)$ is a singleton. Intuitively, $F(k)$ is the object to be deformed, and $F(A)$ is the set of (isomorphism classes of) deformations over $\Spec A$, which can be visualized as a nilpotent thickening of the point $\Spec k$.

Deformation functors often arise as infinitesimal local versions of some moduli functor: given a functor $F\colon\mathsf{Sch}_k^\op\to\mathsf{Set}$ with a distinguished point $p\in F(\Spec k)$, we define the \emph{formal completion} $\widehat{F}$ of $F$ by
\begin{equation*}
\widehat{F}(A):=\{q\in F(\Spec A):q|_{\Spec k}=p\},
\end{equation*}
where $A$ is a local Artin $k$-algebra with maximal ideal $\m_A$ and $q|_{\Spec k}$ denotes the image of $q$ under the map $F(\Spec A)\to F(\Spec A/\m_A)\cong F(\Spec k)$. This is evidently a deformation functor. In particular, if the target category of $F$ is $\mathsf{Ab}$, then the completion of $F$ with respect to the identity is given by
\begin{equation*}
\widehat{F}(A):=\ker(F(\Spec A)\to F(\Spec k)).
\end{equation*}
This is the case with which we shall be primarily concerned.

\begin{example}
\label{ex:completion}
(i) The formal completion of the group scheme $\bG_a=\Spec k[t]$, regarded as its functor of points, is given by
\begin{equation*}
\widehat{\bG}_a(A)=\ker(A\to k)=\m_A.
\end{equation*}
This agrees with our earlier definition in Example~\ref{ex:fg}(i) on the subcategory $\mathsf{Art}_k$; we may therefore regard it as the same formal group.

(ii) Similarly, the formal completion of the group scheme $\bG_m=\Spec k[t,t^{-1}]$ is given by
\begin{equation*}
\widehat{\bG}_m(A)=\ker(A^\times\to k^\times)=(1+\m_A)^\times.
\end{equation*}

(iii) Formally completing an elliptic curve $E/k$ at its identity element yields a formal group $\widehat{E}$ of height $1$ or $2$ (when $k$ is a field of characteristic $p$); in the former case, $E$ is termed \emph{ordinary}, and in the latter, it is termed \emph{supersingular}. Here, $\widehat{E}$ may be viewed as describing the addition law of $E$ in an infinitesimal neighborhood of the identity. After putting $E$ in Weierstrass form (i.e., (\ref{eqn:weierstrass}) with $a_i(t)\in k(t)$ replaced by $a_i\in k$), this intuition can be converted into an algorithm for computing a formal group law of $\widehat{E}$ (see for instance \cite[\S IV.1]{silverman}), which yields
\begin{equation*}
x+_{\widehat{E}}y=x+y+a_1xy-a_2(x^2y+xy^2)+(2a_3x^3y+(a_1a_2-3a_3)x^2y+2a_3xy^3)+\cdots\in\Z[a_1,a_2,a_3,a_4,a_6]\doubles{a,y}.
\end{equation*}
Since $E$ is canonically isomorphic to its Jacobian variety $\Pic^0_E$, the formal group $\widehat{E}$ can be identified with the formal completion $\widehat{\Pic}_E$ of the Picard scheme $\Pic_E$, i.e., the \emph{formal Picard group} of $E$.

(iv) More generally, formally completing any abelian variety at its identity element yields a formal group.
\end{example}

A deformation functor $F$ is \emph{pro-representable} if there exists a local $k$-algebra $R$ with residue field $k\cong R/\m_R$ and finite-dimensional Zariski tangent space $(\m_R/\m_R^2)^\vee$ such that $F$ is represented by $R$ in the category of local $k$-algebras. In this case $F$ is also pro-representable by the $\m_R$-adic completion $\widehat{R}$ of $R$. We would like to determine when a deformation functor $F$ is pro-representable. To understand $F$, one must understand how deformations in $F(A)$ can be lifted across surjections $A'\twoheadrightarrow A$ in $\mathsf{Art}_k$. It suffices to consider \emph{small extensions}, where the kernel $I$ of the quotient $A'\twoheadrightarrow A$ satisfies $I\cdot\m_{A'}=0$: indeed, one easily shows that for any local Artin $k$-algebra $A$, there exists a finite sequence
\begin{equation*}
A=A_0\twoheadrightarrow A_1\twoheadrightarrow\cdots\twoheadrightarrow A_n\twoheadrightarrow A_{n+1}=0.
\end{equation*}
of small extensions (see for instance \cite[Lem.~3.5]{stienstra}). A \emph{tangent-obstruction theory} for $F$ is then defined as the data of two finite-dimensional $k$-vector spaces $T_1$ and $T_2$ such that for any small extension $I\to A'\to A$ in $\mathsf{Art}_k$, there is an exact sequence
\begin{equation}
\label{eqn:tanob}
T_1\otimes_kI\to F(A')\to F(A)\to T_2\otimes_kI,
\end{equation}
functorial in the small extension (see \cite[Rmk. 6.1.20(2)]{fgaexplained}), and whose first map is injective if $A=k$. We may now state Schlessinger's criterion for pro-representability, in the more convenient language of \cite[Cor.~6.3.5]{fgaexplained}:

\begin{theorem}[{\cite[Thm.~2.11]{schlessinger}}]
\label{thm:schlessinger}
A deformation functor $F$ is pro-representable if and only if it admits a tangent-obstruction theory such that the first map of the sequence (\ref{eqn:tanob}) is always injective.
\end{theorem}

Now, we say that a deformation functor $F$ is \emph{formally smooth} if for any surjection $A'\twoheadrightarrow A$ of local Artin $k$-algebras, the map $F(A')\to F(A)$ is surjective. By \cite[Prop.~2.5(i)]{schlessinger}, if $F$ is formally smooth and pro-representable, then it is pro-represented by a power series ring over $k$ in $\dim_k F(k[\epsilon]/\epsilon^2)$ variables. Thus, if $F$ satisfies the conditions of Theorem~\ref{thm:schlessinger}, then it is pro-represented by a power series ring in one variable if and only if $T_2=0$ and $\dim_kT_1=1$ (as the small extension $k[\epsilon]/\epsilon^2\twoheadrightarrow k$ identifies $T_1\cong F(k[\epsilon]/\epsilon^2)$). If in addition $F$ has $\mathsf{Ab}$ as its target category, then $F$ is a ($1$-dimensional) formal group by definition. We will use this criterion in the next section to establish that the formal Brauer group of a K3 surface is indeed a formal group.

\section{The Formal Brauer Group}
\label{sec:4n}

We begin by constructing the formal Brauer group of a K3 surface, following \cite{k3} (see \cite{artinmazur} for the original source). First, the \emph{Brauer group} of a K3 surface $X/k$ is the \'etale cohomology group
\begin{equation*}
\Br(X):=H_\et^2(X,\bG_m)
\end{equation*}
(for an equivalent definition in terms of equivalence classes of Azumaya algebras on $X$, see \cite[\S 18.1.1]{k3}). The \emph{formal Brauer group} $\hBr_X$ of $X$ is then defined as the formal completion of the moduli functor $\Br(X\times -)$, that is,
\begin{equation*}
\hBr_X(A):=\ker(\Br(X\times A)\to \Br(X)).
\end{equation*}
We may likewise define the \emph{formal Picard group} $\hPic_X$ using the Picard group $\Pic(X)=H_\et^1(X,\bG_m)$; this agrees with the definition given in Example~\ref{ex:completion}(iii).

We claim that $\hBr_X$ is a ($1$-dimensional) formal group. Given any small extension $I\to A'\to A$, we have a short exact sequence
\begin{equation*}
0\to(1+\sO_X\otimes_{\underline{k}}\underline{I})^\times\to\pi_*(\bG_m)_{A'}\to\pi_*(\bG_m)_A\to 0
\end{equation*}
of \'etale sheaves on $X$; here $\underline{k}$ and $\underline{I}$ denote the associated constant sheaves, $(\bG_m)_A$ denotes the multiplicative group over $X\times\Spec A$, and $\pi$ denotes the projection to $X$ (and likewise for $A'$). This gives rise to a long exact sequence
\begin{equation*}
\cdots\to\hPic_X(A')\to\hPic_X(A)\to H_\et^2(X,\sO_X)\otimes_kI\to\hBr_X(A')\to\hBr_X(A)\to H_\et^3(X,\sO_X)\otimes_kI\to\cdots.
\end{equation*}
Note that $(1+\sO_X\otimes_{\underline{k}}\underline{I})^\times\cong\sO_X\otimes_{\underline{k}}\underline{I}$ as $I^2=0$. Now, since $H^1(X,\sO_X)=0$, the Picard scheme $\Pic_X$ is $0$-dimensional and smooth (see \cite[Cor.~9.5.13]{fgaexplained}), hence formally smooth, and so $H_\et^2(X,\sO_X)\otimes_kI\to\hBr_X(A')$ is injective. Thus, a tangent-obstruction theory for $\hBr_X$ is given by $T_1:=H_\et^2(X,\sO_X)=H^2(X,\sO_X)$ and $T_2:=H_\et^3(X,\sO_X)=H^3(X,\sO_X)=0$. Since $\dim_kH^2(X,\sO_X)=1$ by Serre duality, it follows from the previous section that $\hBr_X$ is a formal group.

In contrast to the case of the formal Picard group of an elliptic curve (see Example~\ref{ex:completion}(iii)), the formal Brauer group of a K3 surface $X$ may have height $\height(\hBr_X)=1,\ldots,10$ or $\height(\hBr_X)=\infty$; in the latter case, $X$ is termed \emph{supersingular}. Mazur originally proved this in \cite{mazur} by showing that the Newton polygon of the crystalline cohomology $H_\cris^*(X)$ of $X$ lies above the Hodge polygon of $X$; this readily implies that the height of the $F$-isocrystal $H_\cris^2(X)$ is one of these eleven values, and this is known to agree with the height of $\hBr_X$ (see \cite{artinmazur}, or \cite[\S 18.3.3]{k3} for a sketch of the argument). For examples of K3 surfaces realizing each of these heights except $7$, see \cite{yui,goto}.

In the next two sections, we present algorithms of Steinstra \cite{stienstra} and Artin \cite{artin} for computing formal group laws of $\hBr_X$. In his paper, Stienstra identifies $\hBr_X$ with the functor
\begin{equation*}
A\mapsto H_\et^2(X,(1+\sO_X\otimes_{\underline{k}}\underline{\m_A})^\times).
\end{equation*}
This follows as above from the long exact sequence on \'etale cohomology associated to the short exact sequence
\begin{equation*}
0\to(1+\sO_X\otimes_{\underline{k}}\underline{\m_A})^\times\to\pi_*(\bG_m)_A\to\bG_m\to 0.
\end{equation*}
Following Steinstra, we denote the functor $A\mapsto(1+\sF\otimes_{\underline{k}}\underline{\m_A})^\times$ by $\bG_{m,\sF}$ for any sheaf $\sF$ of $k$-algebras on $X$; if in particular $\sF=\sO_X$, then we denote it by $\bG_{m,X}$. On the other hand, Artin regards $\hBr_X$ as the formal completion of the functor $R^2\pi_*\bG_m$ at its zero-section, where $\pi\colon X\to\Spec k$ is the structure map. Since $R^2\pi_*\bG_m$ is the constant sheaf on $\Spec k$ with value $H_\et^2(X,\bG_m)$, this clearly agrees with our previous definition.

\section{Stienstra's Algorithm}
\label{sec:5}

In \cite{stienstra}, Stienstra shows how to compute the logarithm of a formal group arising from a complete intersection in projective space or a double cover of projective space. As we have seen in Examples~\ref{ex:k3}(i)--(iii), K3 surfaces can often be described in this manner. For simplicity, we specialize Steinstra's theorems to the case of the formal Brauer group of a K3 surface (although they hold in much greater generality). Stienstra's methods can be extended to apply to weighted diagonal and quasi-diagonal K3 surfaces; for details, see \cite{yui}.

\begin{theorem}[{\cite[Thm.~1]{stienstra}}]
\label{thm:stienstra}
Let $K$ be a ring that is flat over $\Z$. Let $f_1,\ldots,f_{n-2}$ be a regular sequence of homogeneous polynomials in $K[x_0,\ldots,x_n]$ such that the ideal $(f_1,\ldots,f_{n-2})$ cuts out a K3 surface $X\subset\P_K^n$. Then there is a formal group law for $\hBr_X$ whose logarithm is given by the power series
\begin{equation*}
\log_{\hBr_X}(t)=\sum_{m\ge 1}\frac{1}{m}\beta_mt^m,
\end{equation*}
where $\beta_m$ denotes the coefficient of $(x_0\cdots x_n)^{m-1}$ in $(f_1\cdots f_{n-2})^{m-1}$.
\end{theorem}
\begin{remark}
(i) An implementation of this algorithm can be found in \cite{code} (specifically, see the four functions beginning with ``\texttt{StienstraQuartic}'').

(ii) As in \S\ref{sec:3}, we must work over a $\Z$-flat ring in order to have a logarithm for $\hBr_X$. Although we have previously only worked with K3 surfaces over fields, all previous constructions are easily extended to apply to rings. In general, performing these constructions for a K3 surface $X$ over a scheme $S$ (i.e., a smooth proper morphism $X\to S$ whose geometric fibers are K3 surfaces in the original sense) will be the same as performing them on its fibers $X_k$ in the original sense. For instance, in Example~\ref{ex:fermatsteinstra} below, we compute $\hBr_X$ over $\Z$; reducing this modulo $p$ gives the formal Brauer group of $X_{\Fp}$.
\end{remark}
\begin{proof}
We first claim that the inclusion $\pi\colon X\hookrightarrow\P_K^n$ induces a natural isomorphism
\begin{equation}
\label{eqn:isomone}
\hBr_X=H_\et^2(X,\hbG_{m,X})\cong H_\et^2(\P_K^n,\hbG_{m,\pi_*\sO_X}).
\end{equation}
Since $\pi$ is finite, one has $R^j\pi_*\hbG_{m,X}=0$ for every $j\ge 1$, and therefore the Leray spectral sequence
\begin{equation*}
E_2^{i,j}=H_\et^i(\P_K^n,R^j\pi_*\hbG_{m,X}(A))\implies H_\et^{i+j}(X,\hbG_{m,X}(A))
\end{equation*}
implies that
\begin{equation*}
H_\et^2(X,\hbG_{m,X})\cong H_\et^2(\P_K^n,\pi_*\hbG_{m,X}).
\end{equation*}
Since $\pi$ is affine, the projection formula yields an isomorphism $\pi_*\hbG_{m,X}\cong\hbG_{m,\pi_*\sO_X}$, which proves the claim.

Next, we compute the right-hand side of (\ref{eqn:isomone}) via a kind of Koszul resolution. For any $S\subseteq\{1,\ldots,n-2\}$, let $\sI_S$ be the ideal sheaf on $\P_K^n$ corresponding to $\prod_{i\in S}f_i$. Given two such subsets $S$ and $S'$, define the homomorphism $\partial_{S,S'}\colon\hbG_{m,S}\to\hbG_{m,S'}$ to be $(-1)^k$ times (via the group structure) the homomorphism induced by the inclusion $\sI_S\subset\sI_{S'}$ if $S=\{i_1<i_2<\cdots<i_{\#S}\}$ and $S'=S-\{i_k\}$, and $0$ otherwise. Then we have a complex
\begin{equation}
\label{eqn:koszul}
0\to\bigoplus_{\#S=n-2}\hbG_{m,\sI_S}\to\bigoplus_{\#S=n-3}\hbG_{m,\sI_S}\to\cdots\to\bigoplus_{\#S=1}\hbG_{m,\sI_S}\to\hbG_{m,\P_K^n}\to\hbG_{m,\pi_*\sO_X}\to 0,
\end{equation}
where the final differential is induced by $\pi^\sharp\colon\pi_*\sO_X\to\sO_{\P_K^n}$ and all others are given by the matrices $(\partial_{S,S'})$. We claim that this complex, call it $C^\bullet$, is exact, i.e., $C^\bullet(A)$ is exact for any local Artin $K$-algebra $A$. Since the functor $\sF\mapsto(1+\sF)^\times$ is exact, and $\sO_X$ and all of the sheaves $\sI_S$ are flat over $K$, any small extension $I\to A'\to A$ gives rise to a short exact sequence
\begin{equation*}
0\to C^\bullet(I)\to C^\bullet(A')\to C^\bullet(A)\to 0.
\end{equation*}
The long exact sequence on homology then shows that $C^\bullet(A')$ is exact if $C^\bullet(I)$ and $C^\bullet(A)$ are exact, which reduces us to proving that $C^\bullet(K[\epsilon]/\epsilon^2)$ is exact. Since $\epsilon^2=0$, the group $\hbG_{m,\sI_S}(K[\epsilon]/\epsilon^2)$ is simply $\sI_S\otimes_KK\epsilon$ under addition; thus, $C^\bullet(K[\epsilon]/\epsilon^2)$ is just a sheaf-theoretic version of the Koszul complex
\begin{equation*}
0\to\bigoplus_{\#S=n-2}\sI_S\to\bigoplus_{\#S=n-3}\sI_S\to\cdots\to\bigoplus_{\#S=1}\sI_S\to\sO_{\P_K^n}\to\pi_*\sO_X\to 0
\end{equation*}
associated to $\pi^\sharp$, tensored with $K\epsilon$. The above complex is exact (see for instance \cite[Thm.~7.11]{hartshorne}), and moreover, its terms are flat over $K$, so $C^\bullet(K[\epsilon]/\epsilon^2)$ is exact as well.

Now, observe that $H_\et^i(\P_K^n,\hbG_{m,S})=0$ for every $S\subsetneq\{1,\ldots,n-2\}$ and $i=1,\ldots,n$. Indeed, we may reduce as before by inducting along small extensions to showing that $H^i(\P_K^n,\hbG_{m,S}(K\epsilon))=0$, where $\epsilon^2=0$. Letting $\Ann\epsilon\subset K$ be the annihilator of $\epsilon$, we have an isomorphism $\overline{K}:=K/\Ann\epsilon\cong K\epsilon$ of $K$-modules, and therefore
\begin{equation*}
H_\et^i(\P_K^n,\hbG_{m,S}(K\epsilon))\cong H_\et^i(\P_{\overline{K}}^n,\sO_{\P_{\overline{K}}^n}(-\textstyle\sum_{j\in S}\deg f_j))=0
\end{equation*}
for $i=1,\ldots,n$ as $\sum_{j\in S}\deg f_j\le n$. Thus, we may use (\ref{eqn:koszul}) as an acyclic resolution for computing the cohomology of $\bigoplus_{\#S=n-2}\hbG_{m,\sI_S}$; combining this with (\ref{eqn:isomone}), we obtain a natural isomorphism
\begin{equation*}
H_\et^2(X,\hbG_{m,X})\cong H_\et^n(\P_K^n,\hbG_{m,\sI}),
\end{equation*}
where $\sI:=\sI_{\{1,\ldots,n-2\}}=\sO_{\P_K^n}(-n-1)$ as $\deg f_1\cdots f_n=n+1$ by Example~\ref{ex:k3}(ii). Furthermore, since $H_\et^i(U,\hbG_{m,\sI}(A))=0$ for all $i\ge 1$, affine open $U\subset\P_K^n$, and local Artin $K$-algebras $A$, we may compute this latter group via \v{C}ech cohomology; indeed, arguing as before reduces us to showing this with the quasicoherent sheaf $\sI\otimes_{\underline{K}}\underline{K\epsilon}$ in place of $\hbG_{m,\sI}(A)$.

We now define a coordinatization for $\hBr_X$. For any local Artin $K$-algebra $A$, we may identify $\hbG_{m,\sI}(A)$ with the sheaf $\sI\otimes_{\underline{K}}\underline{\m_A}$ whose group structure is given by the formal group law $x+y-xy$. Letting $\cU=\{U_0,\ldots,U_n\}$ be the standard affine open cover of $\P_K^n$ (that is, $U_i$ is the open subset where the coordinate $t_i$ is invertible), we define a functorial map
\begin{equation}
\label{eqn:coord}
\widehat{\A}_K^1\to\check{H}_\et^n(\cU,\hbG_{m,\sI})
\end{equation}
by sending $a\in\widehat{\A}_K^1(A)=\m_A$ to the \v{C}ech $n$-cocycle given by $f_1\cdots f_{n-2}\cdot t_0^{-1}\cdots t_n^{-1}\otimes a$ on $U_0\cap\cdots\cap U_n$. One easily shows (by reducing as usual to finitely generated Artin $K$-algebras and then inducting along small extensions) that this is a bijection. Similarly, we may identify $\hbG_{m,\sI}(A)$ with the sheaf $\sI\otimes_K\m_A$ whose group structure is given by the formal group law $x+y$; the same map again defines a coordinatization of $\check{H}_\et^n(\cU,\hbG_{a,\sI})$, though now the corresponding formal group law on $\widehat{A}_K^1$ is evidently given by $x+y$.

Finally, we compute the logarithm of $\hBr_X$. Since $K\hookrightarrow K\otimes\Q$ by assumption, we have an isomorphism $\hbG_m\xrightarrow{\cong}\hbG_a$ of formal groups given, in terms of the coordinatizations corresponding to the formal group laws $x+y-xy$ and $x+y$, respectively, by the power series $\log_{\hbG_m}(t)=\sum_{m\ge 1}m^{-1}t^m$ (see Example~\ref{ex:log}). This induces an isomorphism of formal groups over $K\otimes\Q$:
\begin{equation}
\label{eqn:isomma}
\begin{tikzcd}
\check{H}_\et^n(\cU,\hbG_{m,\sI})\arrow{r}{\cong}\arrow{d}{\cong}&\check{H}_\et^n(\cU,\hbG_{a,\sI})\arrow{d}{\cong}\\
\widehat{\A}_{K\otimes\Q}^1\arrow[r,dashed]&\widehat{\A}_{K\otimes\Q}^1.
\end{tikzcd}
\end{equation}
The logarithm of $\hBr_X$, with respect to the coordinatizations of the previous paragraph (i.e., the vertical arrows), is expressed by the dashed arrow. Specifically, it is the power series $\log_{\hBr_X}(t)$ such that the following $n$-cocycles in $A=k\doubles{t}$ are cohomologous:
\begin{gather*}
f_1\cdots f_{n-2}\cdot t_0^{-1}\cdots t_n^{-1}\otimes\log_{\hBr_X}(t),\\
\log_{\hbG_m}(f_1\cdots f_{n-2}\cdot t_0^{-1}\cdots t_n^{-1}\otimes t)=\sum_{m\ge 1}f_1\cdots f_{n-2}\cdot t_0^{-m}\cdots t_n^{-m}\otimes t^m/m.
\end{gather*}
The expression for $\log_{\hBr_X}(t)$ given in the statement of the theorem is now immediate.
\end{proof}

\begin{example}
\label{ex:fermatsteinstra}
By Theorem~\ref{thm:stienstra}, the formal Brauer group of the Fermat quartic $X\subset\P_\Z^3$ (see Example~\ref{ex:fermat}) has logarithm
\begin{equation*}
\log_{\hBr_X}(t)=\sum_{m\ge 0}\frac{(4m)!}{(m!)^4}\frac{t^{4m+1}}{4m+1}.
\end{equation*}
From here, one readily sees that
\begin{equation*}
\height(\hBr_X)=\begin{cases}
1&\text{if }p\equiv 1\bmod 4,\\
\infty&\text{if }p\equiv 3\bmod 4.
\end{cases}
\end{equation*}
Indeed, in the latter case the monomial $t^p$ never appears in the logarithm. In the former, $p$-typicalizing the logarithm (i.e., removing all monomials that are not powers of $p$) gives $\log_{X,p}:=t+\frac{u}{p}t^p+O(p^2)$ for some $u\not\equiv 0\bmod p$; this series has inverse $\exp_{X,p}:=t-\frac{u}{p}t^p+O(p^2)$, so the $p$-series of $\hBr_X$ is
\begin{equation*}
[p]_X=\exp_{X,p}(p\log_{X,p})=pt+ut^p-up^{p-1}t^p+O(t^{p^2})\equiv ut^p+O(t^{p^2})\bmod p.
\end{equation*}
Alternatively, one could compute the formal group law directly as in (\ref{eqn:fglfromlog}): over $\Z$ this yields
\begin{gather*}
\hBr_X(x,y)=x+y-24x^4y-24xy^4-48x^3y^2-48x^2y^3-1944x^8y-6624x^7y^2-14304x^6y^3\\
-20880x^5y^4-20880x^4y^5-14304x^3y^6-6624x^2y^7-1944xy^8+O(11),
\end{gather*}
which in characteristic $3$ has $p$-series $0+O(11)$ and in characteristic $5$ has $p$-series $-x^5+O(11)$.
\end{example}

Next, we state Stienstra's result on double covers of projective space:

\begin{theorem}[{\cite[Thm.~2]{stienstra}}]
\label{thm:stienstradbl}
Let $K$ be a ring that is flat over $\Z$. Let $f$ be a homogeneous degree-$6$ polynomial in $K[x_0,x_1,x_2]$, so that the double cover of $\P_K^2$ given by the equation $w^2=f$ defines a K3 surface $X$ (in a weighted projective space where $w$ is given degree $3$). Then there is a formal group law for $\hBr_X$ whose logarithm is given by the power series
\begin{equation*}
\log_{\hBr_X}(t)=\sum_{m\ge 1}\frac{1}{m}\beta_mt^m,
\end{equation*}
where $\beta_m$ is $0$ if $m$ is even, and the coefficient of $(x_0\cdots x_n)^{m-1}$ in $(f_1\cdots f_{n-2})^{(m-1)/2}$ otherwise.
\end{theorem}
\begin{remark}
An implementation of this algorithm can be found in \cite{code} (specifically, see the four functions beginning with ``\texttt{StienstraDoublePlane}'').
\end{remark}
\begin{proof}
The proof is almost identical to that of the previous theorem, the key difference being that one must replace $f_1\cdots f_{n-2}$ by $w$ in the coordinatization (\ref{eqn:coord}) of $\hBr_X$. For even $m$, the cocycles specified in this coordinatization are all coboundaries, as $w^m=f^{(m-1)/2}\in K[x_0,x_1,x_2]$; for odd $m$, the proof proceeds as before.
\end{proof}

\begin{example}
By Theorem~\ref{thm:stienstradbl}, the formal Brauer group of the double cover $X\to\P_\Z^2$ with equation
\begin{equation*}
w^2=x_0^6+x_1^6+x_2^6
\end{equation*}
has logarithm
\begin{equation*}
\log_{\hBr_X}(t)=\sum_{m\ge 0}\frac{(3m)!}{(m!)^3}\frac{t^{6m+1}}{6m+1}.
\end{equation*}
As before, one can check that
\begin{equation*}
\height(\hBr_X)=\begin{cases}
1&\text{if }p\equiv 1\bmod 6,\\
\infty&\text{otherwise}.
\end{cases}
\end{equation*}
\end{example}

\section{Artin's Algorithm}
\label{sec:6}

In \S 2 of his seminal 1974 paper \cite{artin}, Artin gives an algorithm for computing the formal Brauer group of an elliptic K3 surface, his stated goal being to show that this is possible in spite of the formal Brauer group being ``defined rather abstractly.'' We will of course be using his algorithm in an essential way, to construct K3 spectra from families of elliptic K3 surfaces. Note that in certain cases, \cite{beukers} also allows one to compute the logarithm of the formal Brauer group; however, we do not say more about this result here.

\begin{theorem}[{\cite[\S 2]{artin}}]
\label{thm:artin}
Let $X$ be an K3 surface over a field $k$ that is elliptically fibered with a section, whose minimal Weierstrass model $X'$ has equation
\begin{equation}
\label{eqn:weierstrassmodel}
y^2+a_1(t)xy+a_3(t)y=x^3+a_2(t)x^2+a_4(t)x+a_6(t),
\end{equation}
where the $a_i(t)\in k[t]$ are polynomials of degree at most $2i$. Then a formal group law for $\hBr_X$ may be computed as follows:
\begin{enumerate}
\item Compute the formal group law $+_{\widehat{X}'}$ of $X'$ as an elliptic curve over $\P^1$, i.e., of the formal completion of $X'$ at its zero-section $C_0\subset X'$ with group law given by the usual addition of points on a plane cubic curve (e.g., via Magma's \texttt{FormalGroupLaw()} function).
\item Set
\begin{equation*}
F:=x/t+_{\widehat{X}'}y/t\in k[t,t^{-1}]\doubles{x,y}.
\end{equation*}
\item Let $B_+$ (resp. $B_-$) denote the sum of all monomials in $F$ of degree greater than (resp. less than) $-1$ in $t$, that is, of the form $t^ix^jy^k$ where $i>-1$ (resp. $i<-1$).
\item Set
\begin{equation*}
F:=x/t+_{\widehat{X}'}y/t+_{\widehat{X}'}(-B_+)+_{\widehat{X}'}(-B_-).
\end{equation*}
\item If $F$ is homogeneous of in degree $-1$ in $t$, that is, a sum of monomials of the form $t^{-1}x^jy^k$, then $tF(x,y)$ is a formal group law for $\hBr_X$; otherwise, return to step (3).
\end{enumerate}
\end{theorem}
\begin{remark}
\label{rmk:artin}
(i) An implementation of this algorithm can be found in \cite{code} (specifically, see the two functions beginning with ``\texttt{ArtinElliptic}'').

(ii) Although we have stated the theorem for a elliptic K3 surfaces admitting a section (as we will always be constructing these surfaces directly via their Weierstrass equations), one may always pass from an elliptic K3 surface to one with a section by taking its relative Jacobian fibration (see for instance \cite[\S 11.4]{k3}); by \cite[Lem.~1.8]{artin}, this does not affect the formal Brauer group.

(iii) In practice, one works modulo a certain degree $n$ (usually $p^h+1$, where $p:=\chr(k)$ and $h$ is the height of interest). Heuristically, around $n/2$ iterations of steps (3)--(5) are usually necessary in this case.

(iv) It is easy to misinterpret the language in \cite[\S 2]{artin} to mean that after step (2) one should simply take the terms in $F$ of degree $-1$ in $t$ and skip to step (5); however, this yields power series that are \emph{not} formal group laws (failing to satisfy associativity). This is because the manual elimination of the monomials not of the form $t^{-1}x^jy^k$ in step (4) can also change those of that form in higher degrees via the non-linear terms of $+_{\widehat{X}'}$.
\end{remark}
\begin{proof}
After labeling morphisms by $X\xrightarrow{f}\P_k^1\xrightarrow{g}\Spec k$, the Leray spectral sequence gives
\begin{equation*}
E_2^{i,j}=R^ig_*R^jf_*\bG_m\implies R^{i+j}(g\circ f)_*\bG_m,
\end{equation*}
so the terms which may contribute to $R^2(g\circ f)_*\bG_m$ are $E_2^{2,0}$, $E_2^{1,1}$, and $E_2^{0,2}$. Since $f_*\bG_m=\bG_m$ and $\P_k^1$ has dimension $1$, the fundamental exact sequence implies that $E_2^{2,0}=R^2g_*\bG_m$ is discrete. Likewise, $E_2^{0,2}=g_*R^2f_*\bG_m$ is discrete as $f$ has relative dimension $1$. Thus, the formal structure of $R^2(g\circ f)_*\bG_m$ at its zero-section is that of $E_2^{1,1}=R^1g_*R^1f_*\bG_m=R^1g_*\Pic_{X/\P_k^1}$, which is in turn the same as that of $R^1g_*X$, as one can show. Thus, $\hBr_X$ is the formal completion of the functor $R^1g_*X$ at its zero-section, hence by \cite[Prop.~II.1.7]{artinmazur}, there is an isomorphism
\begin{equation*}
\hBr_X(A)\cong H_\et^1(\P_k^1,\pi_*\widehat{X}_A),
\end{equation*}
natural in the local Artin $k$-algebra $A$. Here $\widehat{X}$ denotes the formal completion of $X$ at its zero-section $C_0$, $\widehat{X}_A:=\widehat{X}\times\Spec A$, and $\pi\colon(\P_k^1)_A\to\P_k^1$ denotes projection onto the first factor, which is an isomorphism of underlying topological spaces.

Let $\sN_{C_0/X}$ denote the normal sheaf of $X$ at its zero-section $C_0$. By \cite[Lem.~2.2.1]{k3}, $C_0$ has self-intersection number $(C_0.C_0)=-2$, hence $\sN_{C_0/X}\cong\sO_X(C_0)|_{C_0}\cong\sO_{\P^1}(-2)$. Let $A'\to A$ be a small extension of local Artin $k$-algebras. Then $\widehat{X}_A$ and $\widehat{X}_{A'}$ corepresent sheaves on both the Zariski and \'etale sites of $(\P_k^1)_A$, that fit into a short exact sequence
\begin{equation*}
0\to\sN_{C_0/X}\to\pi_*\widehat{X}_{A'}\to\pi_*\widehat{X}_A\to 0
\end{equation*}
of sheaves for both the Zariski and \'etale topologies on $\P_k^1$. Intuitively, sections of $\sN_{C_0/X}$ correspond to first-order deformations of $C_0$ inside $X$, and sections of $\pi_*\widehat{X}_A$ correspond to $\len(A)$th-order deformations of $C_0$ inside $X$, where $\len(A)$ denotes the length of $A$ (which here coincides with its dimension over $k$); exactness then corresponds to the fact that $\len(A')=\len(A)+1$, and this is not difficult to verify by working locally on $\P_k^1$. Now, consider the commutative diagram
\begin{equation}
\label{eqn:zaretcomp}
\begin{tikzcd}[column sep=1em]
0\arrow{r}&H^0(\pi_*\widehat{X}_{A'})\arrow{r}\arrow{d}&H^0(\pi_*\widehat{X}_A)\arrow{r}\arrow{d}&H^1(\sN_{C_0/X})\arrow{r}\arrow{d}&H^1(\pi_*\widehat{X}_{A'})\arrow{r}\arrow{d}&H^1(\pi_*\widehat{X}_A)\arrow{r}\arrow{d}&0\\
0\arrow{r}&H_\et^0(\pi_*\widehat{X}_{A'})\arrow{r}&H_\et^0(\pi_*\widehat{X}_A)\arrow{r}&H_\et^1(\sN_{C_0/X})\arrow{r}&H_\et^1(\pi_*\widehat{X}_{A'})\arrow{r}&H_\et^1(\pi_*\widehat{X}_A)\arrow{r}&0
\end{tikzcd}
\end{equation}
with exact rows, where all cohomology is taken over $\P_k^1$. Indeed, since $\sN_{C_0/X}$ is quasi-coherent and of negative degree, we have $H_\et^i(\P_k^1,\sN_{C_0/X})=H^i(\P_k^1,\sN_{C_0/X})=0$ for $i=0,2$. Since $\pi_*\widehat{X}_k$ is the zero sheaf (intuitively, the only ``zeroth-order'' deformation of a Zariski or \'etale open of $\P_k^1$ is the trivial one), inducting along small extensions in the three left-most columns of (\ref{eqn:zaretcomp}) shows that $H_\et^0(\pi_*\widehat{X}_A)=H^0(\pi_*\widehat{X}_A)=0$ for all local Artin $k$-algebras $A$. Applying an analogous argument to the five right-most columns of (\ref{eqn:zaretcomp}) then shows that
\begin{equation*}
H_\et^1(\P_k^1,\pi_*\widehat{X}_A)\cong H^1(\P_k^1,\pi_*\widehat{X}_A)
\end{equation*}
for any local Artin $k$-algebra $A$, which may be computed as \v{C}ech cohomology for any affine covering of $\P_k^1$.

To this end, choose homogeneous coordinates $t_0$ and $t_1$ on $\P_k^1$, and let $\cU=\{U_0,U_1\}$ be the standard affine open cover thereof (i.e., $U_i$ is the open subset where $t_i$ is invertible). We first compute the \v{C}ech cohomology group $\check{H}^1(\cU,\sN_{C_0/X})$: since $\cU$ is a trivializing open cover for the invertible sheaf $\sN_{C_0/X}\cong\sO_{\P^1}(-2)$, we may choose bases $\{u_0\}$ and $\{u_1\}$ for $\sN_{C_0/X}$ over $U_0$ and $U_1$, respectively, related over $U_0\cap U_1=\Spec k[t_{1/0},t_{0/1}]$ by $u_0=t_{1/0}^2u_1$, where we have let $t_{i/j}:=t_i/t_j$. A $1$-cocycle for $\cU$ is any section of $\sN_{C_0/X}$ over $U_0\cap U_1$, say $f(t_{1/0},t_{0/1})u_0$. The $1$-coboundaries are the sections of $\sN_{C_0/X}$ over $U_0\cap U_1$ of the form
\begin{equation}
\label{eqn:cobound}
g_0(t_{1/0})u_0+g_1(t_{0/1})u_1=(g_0(t_{1/0})+t_{0/1}^2g_1(t_{0/1}))u_0,
\end{equation}
by which we may eliminate all terms of the $1$-cocycle $f(t_{1/0},t_{0/1})u_0$ except for the monomial $t_{0/1}u_0$ and its coefficient. Thus, $\check{H}^1(\cU,\sN_{C_0/X})$ is the $1$-dimensional $k$-vector space spanned by $t_{0/1}u_0=t_{1/0}u_1$. Let $\mathbf{V}(\sN_{C_0/X})$ denote the total space of $\sN_{C_0/X}$, which is given by gluing $\A_{U_0}^1=\Spec k[t_{1/0},u_0]$ and $\A_{U_1}^1=\Spec k[t_{0/1},u_1]$ together over $U_0\cap U_1$ via $u_0=t_{1/0}^1u_1$, as before. Its sheaf of sections is isomorphic to $\sN_{C_0/X}$, so that the cohomology class of $t_{0/1}u_0=t_{1/0}u_1$ corresponds to that represented by the $1$-cocycle $\Spec k[t_{1/0},t_{0/1}]\to\mathbf{V}(\sN_{C_0/X})$ given by $u_0\mapsto t_{0/1}$ and $u_1\mapsto t_{1/0}$ over $U_0\cap U_1$.

Similarly, for the \v{C}ech cohomology group $\hBr_X(A)\cong\check{H}^1(\cU,\pi_*\widehat{X}_A)$, a $1$-cocycle is given by a section $(U_0\cap U_1)_A\to X_A$ deforming the zero-section over $U_0\cap U_1$, by the definition of formal completion. Since $\widehat{X}$ is a $1$-dimensional formal group over $\P_k^1$, this is equivalent, after a choice of coordinate, to a map
\begin{equation*}
k[t_{1/0},t_{0/1}]\doubles{x}\to k[t_{1/0},t_{0/1}]\otimes A=A[t_{1,0},t_{0,1}]
\end{equation*}
of local Artin $k[t_{1/0},t_{0/1}]$-algebras (i.e., such that the image of $x$ vanishes modulo $\m_A$). Now, $\hBr_X$ is a $1$-dimensional formal group, so it must admit a coordinatization, that is, an isomorphism
\begin{equation*}
\check{H}^1(\cU,\pi_*\widehat{X}_A)\cong\Hom_{\mathsf{Art}_k}(k\doubles{x},A)
\end{equation*}
natural in $A$; this furnishes a universal cohomology class in $\check{H}^1(\cU,\pi_*\widehat{X}_{k\doubles{x}})$ corresponding to $1_{k\doubles{x}}$ on the right-hand side. Since $\sN_{C_0/X}$ is the first-order approximation to $\pi_*\widehat{X}_A$, this must be represented by the $1$-cocycle
\begin{align*}
k[t_{1/0},t_{0/1}]\doubles{x}&\to k[t_{1/0},t_{0/1}]\doubles{x},\\
x&\mapsto t_{0/1}x.
\end{align*}
Indeed, this must be true modulo $(x^2)$; any monomials in the image of $x$ whose degree in $t_{0/1}$ is not $1$ may be eliminated via coboundaries as in (\ref{eqn:cobound}); and a change of coordinates in $k\doubles{x}$ brings the remaining power series to $t_{0/1}x$, as desired.

Summing this $1$-cocycle with itself in $\pi_*\widehat{X}_{k[x_1,x_2]}$ gives the $1$-cocycle
\begin{equation}
\label{eqn:sum}
\begin{split}
k[t_{1/0},t_{0/1}]\doubles{x}&\to k[t_{1/0},t_{0/1}]\doubles{x_1,x_2},\\
x&\mapsto t_{0/1}x_1+_{\widehat{X}}t_{0/1}x_2,
\end{split}
\end{equation}
where $+_{\widehat{X}}$ denotes the formal group law of $\widehat{X}$. Since $X$ and $X'$ have the same formal structure at their zero-sections, we may compute this as in step (1) in the statement of the theorem. By universality, the $1$-cocycle (\ref{eqn:sum}) corresponds to a map $k\doubles{x}\to k\doubles{x_1,x_2}$, that is, a power series $\hBr_X(x_1,x_2)\in k\doubles{x_1,x_2}$, such that the $1$-cocycle $x\mapsto t_{0/1}\hBr_X(x_1,x_2)$ represents the same cohomology class as (\ref{eqn:sum}) in $\check{H}^1(\cU,\pi_*\widehat{X}_{k\doubles{x_1,x_2}})$. This power series is by definition a formal group law for $\hBr_X$. As before, the $1$-coboundaries in $\pi_*\widehat{X}_{k\doubles{x_1,x_2}}$ are sections over $U_0\cap U_1$ of the form $x\mapsto g_0+_{\widehat{X}}t_{0/1}^2g_1$, where $g_0\in k[t_{1/0}]\doubles{x_1,x_2}$ and $g_1\in k[t_{0/1}]\doubles{x_1,x_2}$ are power series such that $g_0(0,0)=0$ and $g_1(0,0)=0$. Since the formal group law $+_{\widehat{X}}$ is ordinary addition modulo higher-order terms, we may inductively eliminate all monomials of $t_{0/1}x_1+_{\widehat{X}}t_{0/1}x_2$ that are not of the form $t_{0/1}x_1^ix_2^j$ using coboundaries; one relatively efficient way of doing so is described in steps (3)--(5) in the statement of the theorem (note that there we have replaced $t_{1/0}$, $x_1$, and $x_2$ by $t$, $x$, and $y$, respectively). Multiplying the resultant power series by $t_{1/0}$ leaves $\hBr_X(x_1,x_2)\in k\doubles{x_1,x_2}$, as desired.
\end{proof}

\begin{example}
Applying Theorem~\ref{thm:artin} to the minimal Weierstrass model (\ref{eqn:wmfq}) for the Fermat quartic $X$ in characteristic $5$, we find that
\begin{align*}
\hBr_X(x,y)&=x+y+2x^2y+2xy^2+4x^3y^2+4x^2y^3+x^6y+3x^5y^2+3x^4y^3+3x^3y^4+2x^2y^5+xy^6\\
&+x^8y+2x^7y^2+3x^6y^3+3x^3y^6+2x^2y^7+xy^8+O(11),
\end{align*}
which has $p$-series $4x^5+O(11)$. Thus, $\height(\hBr_X)=1$, which agrees with Example~\ref{ex:fermatsteinstra}.
\end{example}

Using Theorem~\ref{thm:artin}, we can obtain a general characterization of the height of an elliptic K3 surface in terms of the coefficients of its minimal Weierstrass model. Indeed, simply compute the formal group law over $k[\{a_{ij}\}]$ corresponding to the Weierstrass equation given by $a_i(t)=\sum_{j=0}^{2i}a_{i,j}t^j$ for each $i=1,2,3,4,6$; then compute the coefficient of $x^{p^h}$ in its $p$-series for each $h\ge 1$. Naturally, this is much easier for smaller primes $p$ and heights $h$. We note the following correction of \cite[Thm.~2.12]{artin} on the case $p=2$:

\begin{theorem}
In characteristic $2$, the equation (\ref{eqn:weierstrassmodel}) can be chosen so that $a_2(t)=0$ and $a_{1,2}=1$, after a suitable change of coordinates. Then the height $h$ of the elliptic K3 surface it describes satisfies
\begin{itemize}
\item $h=1$ if $a_{1,1}\ne 0$;
\item $h\ge 2$ if $a_{1,1}=0$;
\item $h\ge 3$ if $a_{1,1}=a_{3,3}=0$;
\item $h\ge 4$ if $a_{1,1}=a_{3,3}=0$ and
\begin{gather*}
a_{1,0}^3a_{4,7}+a_{1,0}^2a_{4,5}+a_{1,0}a_{3,1}a_{3,6}+a_{1,0}a_{3,2}a_{3,5}+a_{1,0}a_{4,3}+a_{1,0}a_{6,7}+a_{3,0}a_{3,5}+a_{3,0}a_{4,7}\\
+a_{3,1}a_{3,4}+a_{3,1}a_{4,6}+a_{3,2}a_{4,5}+a_{3,4}a_{4,3}+a_{3,5}a_{4,2}+a_{3,6}a_{4,1}+a_{4,1}+a_{6,5}=0.
\end{gather*}
\end{itemize}
\end{theorem}
In Artin's original statement of the theorem, the condition for $h\ge 4$ is given instead by $a_{1,1}=a_{3,3}=a_{6,5}=0$, but this is incorrect (despite being a nicer statement). For example, the K3 surface in characteristic $2$ with Weierstrass equation
\begin{equation*}
y^2+t^2xy=x^3+tx
\end{equation*}
has formal group law $x+y+x^4y^4+O(9)$ and $p$-series $x^8+O(9)$, hence has height $3$, rather than $4$ as predicted by Artin's original statement.

\section{The Landweber Exact Functor Theorem}
\label{sec:7}

Having explained how to compute the formal Brauer group of a K3 surface, we now turn to the problem of using these to construct homology theories, or equivalently, spectra. Let $\bG$ be a formal group over a ring $R$. Choose a coordinate $t$ on $\bG$ and let
\begin{equation*}
[p]_{\bG}=a_1t+a_2t^2+\cdots+a_pt^p+a_{p+1}t^{p+1}+\cdots
\end{equation*}
be its $p$-series. For each integer $n\ge 0$, let $v_n:=a_{p^n}$; thus, $v_0=a_1=p$. Although this definition depends on the choice of coordinate $t$, the ideals $(p,v_1,\ldots,v_n)$ do not (see for instance \cite[\S21]{coctalos}). Thus, we may define the formal group $\bG$ to be \emph{regular at $p$} if the sequence $p,v_1,v_2,\ldots\in R$ is regular, i.e., the map
\begin{equation*}
v_n\colon R/(p,v_1,\ldots,v_{n-1})R\to R/(p,v_1,\ldots,v_{n-1})R
\end{equation*}
is injective for all $n\ge 0$. Replacing $R$ by the graded ring $R[u^{\pm 1}]$ if necessary, we may assume that the formal group law $\bG$ is graded in even degrees, hence classified by a graded morphism $\MU^*\to R$. The Landweber exact functor theorem then states:
\begin{theorem}
\label{thm:landweber}
If $\bG$ is $p$-regular for all primes $p$, then the functor
\begin{equation*}
X\mapsto\MU^*X\otimes_{\MU^*}R
\end{equation*}
is a cohomology theory. It is representable by an even periodic ring spectrum $E$ satisfying $\pi_*E\cong R$ and $\bG_E\cong\bG$.
\end{theorem}
Such a formal group $\bG$ is called \emph{Landweber exact}. Note that it suffices to check Landweber exactness before grading $R$ and $\bG$. This provides a purely algebraic---and often very easily computable---criterion for a formal group to define a cohomology theory. The exact functor theorem also admits an interpretation in the language of stacks: $\bG$ is classified by a map $\Spec R\to\sM_{\mathrm{FG}}$ to the moduli stack of formal groups, which is flat if and only if $\bG$ is Landweber exact. For proofs in this language, see \cite{coctalos,luriechrom}; for a more elementary proof, see Haynes Miller's excellent expository account \cite{haynes}.

In particular, if $\bG$ is the formal Brauer group of a K3 surface over $R$, then the spectrum representing the associated cohomology theory will be a K3 spectrum (see Definition~\ref{def:k3spectrum}). We will use this observation in the next section, but before going on, let us give some examples showing how to check Landweber exactness:
\begin{example}
\label{ex:landweber}
(i) Every formal group law over a $\Q$-algebra $R$ is Landweber exact, as $v_0=p$ is a unit in $R$, and therefore $R/v_0R=0$ (after which the regularity condition becomes trivial). Likewise, any formal group law over a $\Z_{(p)}$-algebra (or a $\Z_p$-algebra) is automatically $q$-regular at all primes $q\ne p$.

(ii) The additive formal group $\widehat{\bG}_a$ is not Landweber exact over $\Z$, as for any prime $p$ we have $v_0=p$ and $v_1=0$, which does not give an injective endomorphism of $\Z/p\Z$. However, it is Landweber exact over $\Q$ by (i), corresponding to the Eilenberg--Maclane spectrum $\H\Q$, i.e., ordinary homology with rational coefficients.

(iii) Consider the multiplicative formal group $\widehat{\bG}_m$ over $\Z$. As in Example~\ref{ex:fg}(ii), it has formal group $x+y-xy$, which we may replace by the graded formal group law $x+y-\beta xy$ over $\Z[\beta^{\pm 1}]$, which determines a map $\MU^*\to\Z[\beta^{\pm 1}]$ of graded rings (here we use $\beta$ instead of $u$, as it stands for ``Bott periodicity''). For any prime $p$, we see that $v_0=p$ is a non-zero-divisor on $\Z[\beta^{\pm 1}]$, and that modulo $p$, the $p$-series of $\widehat{\bG}_m$ is given by $\pm\beta^{p-1}t^p$, so that $v_1\equiv\pm\beta^{p-1}\bmod p$, which acts invertibly on $\Z[\beta^{\pm 1}]/p$. Thus, $\Z[\beta^{\pm 1}]/(v_0,v_1)=0$, so that $\widehat{\bG}_m$ is $p$-regular and therefore Landweber exact. As it turns out, the homology theory $E$ with coefficient ring $\pi_*E\cong\Z[\beta^{\pm 1}]$ given by the Landweber exact functor theorem is complex $K$-theory (see \cite{luriechrom}).

(iv) Examples bearing more similarity to the constructions in the next section are given by elliptic cohomology theories; here we treat the universal case. Let $A:=\Z[a_1,a_2,a_3,a_4,a_6]$ be the ring over which the universal Weierstrass curve lives; after inverting the discriminant $\Delta$, we may consider the universal smooth Weierstrass curve $C$ over the ring $A[\Delta^{-1}]$. We claim that its formal group, in the sense of Example~\ref{ex:completion}(iii), is Landweber exact. Let $p$ be a prime. Since $A[\Delta^{-1}]$ is torsion-free, multiplication by $v_0=p$ is an injective endomorphism of $A[\Delta^{-1}]$. Since $A[\Delta^{-1}]/p$ is an integral domain, multiplication by $v_1$ is injective if and only if the image of $v_1$ in $A[\Delta^{-1}]/p$ is nonzero. Equivalently, the reduction of $C$ modulo $p$ has a fiber that is ordinary; this is certainly true, as most elliptic curves are not supersingular. Finally, the image of $v_2$ in $A[\Delta^{-1}]/(p,v_1)$ must be a unit: otherwise, there exists a maximal ideal $\m\subset A$ containing $p$, $v_1$, and $v_2$, so that the reduction of $C$ defines an elliptic curve over the field $A[\Delta^{-1}]/\m$ with formal group law of height greater than $2$, a contradiction. Thus, the $p$-regularity condition for $v_3,v_4,\ldots$ is trivial, and the formal group of $C$ is Landweber exact. One can use this fact to show that the map $\sM_{\mathrm{Ell}}\to\sM_{\mathrm{FG}}$ that associates to any elliptic curve its formal group is flat (see \cite{tmf}). Here $\sM_{\mathrm{Ell}}$ denotes the moduli stack of elliptic curves, which is associated to the Hopf algebroid $(A[\Delta^{-1}],\Gamma[\Delta^{-1})$, where the ring $\Gamma:=A[u^{\pm 1},r,s,t]$ parametrizes isomorphisms of Weierstrass curves. Altogether, one obtains an elliptic cohomology theory $\mathrm{Ell}$ with coefficient ring $\pi_*\mathrm{Ell}\cong A[\Delta^{-1}]$, where each $a_i$ is graded in degree $2i$.
\end{example}

\section{Explicit K3 Spectra}
\label{sec:8}

We now provide three explicit examples of height-$3$ K3 spectra, each arising from a family of K3 surfaces of a different type discussed in \S\ref{sec:2} to which either Stienstra's or Artin's algorithm applies.

\begin{example}
Our first K3 spectrum $Q$ arises from the following family over $\Spec\Z_3[a,b]$ of quartic surfaces in $\P^3$, with respect to homogeneous coordinates $x_0$, $x_1$, $x_2$, and $x_3$:
\begin{equation}
\label{eqn:famone}
x_0^4 + x_0^2x_1x_3 + x_0x_1x_2^2 + x_0x_3^3 + x_1^4 + x_2^4 + ax_1x_3^3 + bx_1x_2^2x_3=0.
\end{equation}
These quartics are not all smooth, but we can compute the locus in $\Spec\Z_3[a,b]$ where this family is singular via the multivariate resultant of (\ref{eqn:famone}) and its partial derivaties (e.g., using Macaulay2's \texttt{Discriminant()} function \cite{M2}). This is a closed subscheme cut out by a polynomial $\Delta$ in $a$ and $b$; the restriction of this family to $\Spec\Z_3[a,b,\Delta^{-1}]$ is then a family of K3 surfaces by Example~\ref{ex:k3}(i).

Computing the $3$-series of its formal Brauer group via Theorem~\ref{thm:stienstra} and reducing modulo $3$ gives $v_1=-b$ and $v_2=-a^2-ab^2$; the expression for $v_3$ is more complicated, but one can check that it is equivalent to $1$ modulo $(3,v_1,v_2)$ (in particular, this family has only one height-$3$ point modulo $3$, given by $a=b=0$). Thus, the formal Brauer group of this family is $3$-regular over $\Z_3[a,b,\Delta^{-1}]$, hence by Examples~\ref{ex:landweber}(i), it is Landweber exact. Note that one must also check that the closed subschemes of the base scheme cut out by the ideals $(3)$, $(3,v_1)$, and $(3,v_1,v_2)$ are not contained in the singular locus; this can be done by computing $\Delta$ directly, or by verifying that specific points of these loci are non-singular (e.g., via Magma's \texttt{IsSingular()} function).

Applying Theorem~\ref{thm:landweber} then produces a K3 spectrum $Q$ with $\pi_*Q\cong\Z_3[a,b,\Delta^{-1},u^{\pm 1}]$, with $a$ and $b$ graded in degree $0$ and $u$ in degree $2$. Since, as noted above, the height-$3$ locus consists of only one point modulo $3$, we have $K(3)_*Q=K(3)_*=\F_3[v_3,v_3^{-1}]$ with $v_3$ graded in degree $2(3^3-1)=52$, where $K(n)$ denotes Morava $K$-theory.
\end{example}
\begin{remark}
(i) The computations in this section may be verified using the three functions in \cite{code} beginning with ``\texttt{VerifyK3Spectrum}.''

(ii) We have chosen to work locally at the prime $3$ in constructing these spectra for the most part arbitrarily. Computations up to height $3$ at primes greater than $3$ tend to be difficult, but one could of course make analogous constructions at the prime $2$ without any trouble.

(iii) One could also look at the space of \emph{all} smooth quartics in $\P^3$, rather than a specific family centered around a single height-$3$ surface. The smooth quartics in $\P^3$ are parametrized by $\mathbb{P}^{{7\choose 5}-1}=\mathbb{P}^{34}$, and a standard argument via the Veronese embedding shows that removing the discriminant locus $\Delta$ gives an affine scheme $\mathbb{P}^{34}\setminus\Delta$. Taking a stack-quotient by the action of the affine group $\mathrm{PGL}(4)$ via change of coordinates gives a presentation of the $19$-dimensional moduli stack $\sM_{\mathrm{SQ}}$ of smooth quartics as a Hopf algebroid. Applying Theorem~\ref{thm:artin} to the universal smooth quartic over $\sM_{\mathrm{SQ}}$ shows that the sequence $3,v_1,v_2,v_3$ is regular. It is difficult to compute $v_i$ for any $i>3$; however, we may simply invert $v_3$ in $\sM_{\mathrm{SQ}}$ to restrict to the locus of smooth quartics of height at most $3$, which, after localizing at the prime $3$, forces Landweber exactness. There are good reasons not to work in this more general setting, though: computations are generally quite expensive, and the coefficient ring of the resulting spectrum is essentially intractable. Thus, henceforth, we limit ourselves to K3 spectra arising from $2$-dimensional families of K3 surfaces containing a single height-$3$ surface modulo $3$.

(iv) As a result, these three spectra have very similar coefficient rings; however, they are distinguishable by their formal groups and singular loci.
\end{remark}
\begin{example}
Our second K3 spectrum $D$ arises from the following family over $\Spec\Z_3[a,b]$ of double covers of $\P^2$ branched along degree-$6$ plane curves, with respect to homogeneous coordinates $x_0$, $x_1$, and $x_2$ of degree $1$ and $w$ of degree $3$:
\begin{equation*}
w^2=-x_0^6+x_0^2x_1^4+x_0x_1^5+x_1x_2^5+x_2^6+ax_0x_1^2x_2^3+bx_0^2x_1^2x_2^2.
\end{equation*}
By Example~\ref{ex:k3}(iii), it suffices in characteristic $3$ for the plane curves cut out by the right-hand side to be smooth, and we may then proceed to identify the singular locus $\Delta$ as before (note that to actually identify these curves, we have used \cite[Tab.~1]{poonen}). Using Theorem~\ref{thm:stienstradbl}, we obtain $v_1=b$ and $v_2=a$ modulo $3$, as well as $v_3=1$ modulo $(3,v_1,v_2)$, from which Landweber exactness is immediate (and smoothness can be checked as before). Theorem~\ref{thm:landweber} then produces a K3 spectrum $D$ with $\pi_*D\cong\Z_3[a,b,\Delta^{-1},u^{\pm 1}]$, with $|a|=|b|=0$ and $|u|=2$ as before. Since the height-$3$ locus again consists of only one point modulo $3$ (namely $a=b=0$), we have $K(3)_*D=K(3)_*$ as before.
\end{example}

\begin{example}
Our third K3 spectrum $E$ arises from the following family over $\Spec\Z_3[a,b]$ of elliptic K3 surfaces in minimal Weierstrass form (see (\ref{eqn:weierstrass})):
\begin{equation*}
y^2+(a+bt)xy+t^2y=x^3+(1+t)x^2+(1+t^4+t^8)x+(t^7+t^8).
\end{equation*}
Indeed, one can check that its discriminant has $t$-adic valuation at most $3$ as $a$ and $b$ vary (in particular, it never vanishes, so there is no singular locus to exclude). Using Theorem~\ref{thm:artin}, we compute modulo $3$ that $v_1=b^2$ and $v_2=a^4-ab+b^4$. Now, $v_3$ is harder to compute modulo $3$, though modulo $(3,a,b)$ it is equal to $-1$. One can then show that this implies that $v_3$ is a non-zero-divisor in $\Z_3[a,b]/(3,v_1,v_2)$, which verifies Landweber exactness. Finally, Theorem~\ref{thm:landweber} produces a K3 spectrum $E$ with $\pi_*E\cong\Z_3[a,b,u^{\pm 1}]$, with $|a|=|b|=0$ and $|u|=2$ as before. Since the height-$3$ locus again consists of only one point modulo $3$ (namely $a=b=0$), we have $K(3)_*E=K(3)_*$.
\end{example}

\section{Further Questions}
\label{sec:9}

The results in the previous section raise many more questions than they answer; we list some of these questions here.

\subsection{Properties of K3 Spectra}
\label{sec:prop}
We would like to deduce more topological information about the spectra constructed in \S\ref{sec:8}. First, are they $A_\infty$- or $E_\infty$-ring spectra? This could conceivably be answered using either Goerss--Hopkins or Robinson obstruction theory; indeed, the former was used in \cite{szymik2} to show that any K3 spectrum refining a local ring of the moduli stack of \emph{ordinary} $p$-primitively polarized K3 surfaces in characteristic $p$ admits an essentially unique $E_\infty$ structure (these methods do not readily generalize to higher height, for reasons explained in \S\ref{sec:tmf}). If these spectra are indeed $E_\infty$, then what are their power operations, and how do they relate to morphisms of K3 surfaces (for Morava $E$-theory and the spectrum $\mathrm{TMF}$ of topological modular forms, these arise from isogenies of formal groups and elliptic curves, respectively)? What are their localizations at the Morava $K$-theory spectra? These questions are essential to the construction of $\mathrm{TMF}$, and so could have applications to \S\ref{sec:tmf} below.

\subsection{Wild Automorphisms of Morava $E$-Theories} Can one construct a height-$3$ (or greater) K3 spectrum with a $p$-torsion automorphism at an odd prime $p$, or a $4$-torsion automorphism at $2$? This amounts to finding a family of K3 surfaces in characteristic $p$ (resp. $2$) whose height-$3$ locus has a $p$-torsion (resp. $4$-torsion) automorphism (note that in \cite{szymik2}, it was shown that the induced automorphism of the K3 spectrum can be taken to be $E_\infty$ in the ordinary $p$-primitively polarized case, see \S\ref{sec:prop}). This can only occur if $p\le 11$ \cite{dolgachev}; such automorphisms can arise, for instance, via translation by a $p$-torsion section of an elliptic K3 surface (such surfaces do not exist for $p=7,11$, but for $p=3,5$ they do and are generically ordinary, see \cite{ito}). However, efficiency limitations in Artin's algorithm \cite{artin} for computing formal Brauer groups make this difficult. A combination of Artin's algorithm and the algorithm presented by Sch\"utt in \cite{schuett} for computing the height of a K3 surface with a wild automorphism via the characteristic polynomial of the trace of Frobenius on \'etale cohomology, suitably adapted, could yield results. This might allow for extensions of recent work of Hahn and Shi \cite{hahnshi} on orientations of Morava $E$-theories at the prime $2$ compatible with $2$-torsion automorphisms.

\subsection{Topological Modular Forms for K3 Surfaces}
\label{sec:tmf}

We would like to adapt the Lurie--Goerss--Hopkins--Miller theorem \cite{lurie} to the situation of K3 surfaces. Roughly, it states that a formally \'etale map from an algebraic stack $\sX$ to the moduli $\sM_p(h)$ of $1$-dimensional $p$-divisible groups of height $h$ gives rise to a sheaf of even weakly-periodic $E_\infty$-ring spectra on the \'etale site of $\sX$ (for details, see \cite[Thm.~5.1]{lawson}). Let $\sM_{\mathrm{Ell}}$ be the moduli stack of elliptic curves. For elliptic cohomology, Serre--Tate theory guarantees that the map $\sM_{\mathrm{Ell}}\to\sM_p(2)$ sending an elliptic curve to its associated $p$-divisible group is formally \'etale. Thus, one obtains a sheaf of $E_\infty$-ring spectra on $\sM_{\mathrm{Ell}}$, with global sections the famed spectrum $\mathrm{TMF}$ of topological modular forms \cite{tmf}. However, for K3 cohomology, Serre--Tate theory does not apply (as K3 surfaces are not abelian varieties), so some substitute is required (in fact, there is a generalization of Serre--Tate theory to K3 surfaces due to Deligne--Illusie \cite{deligne,deligne2,katz}, but it utilizes $F$-crystals rather than $p$-divisible groups to parameterize the deformations of a K3 surface).

One idea is to use some moduli of Dieudonn\'e modules, rather than $p$-divisible groups, by sending a K3 surface to the Dieudonn\'e module associated to its formal Brauer group. However, the map to this moduli from a chosen Landweber exact moduli $\sX$ of K3 surfaces is not in general formally \'etale, as the formal Brauer group cannot have constant height across $\sX$, and therefore the associated Dieudonn\'e module cannot have constant dimension over the Witt ring $W$ of the base field. Nonetheless, Illusie's slope spectral sequence \cite{illusie} provides an embedding of the Dieudonn\'e module associated to the formal Brauer group of a K3 surface $X$ into its crystalline cohomology group $H_\cris^2(X/W)$, which has constant dimension $22$ and carries a Dieudonn\'e module structure of its own.

Another idea is to impose more structure on the K3 surfaces in question to force their moduli to be formally \'etale over the moduli of $p$-divisible groups. This was accomplished in \cite{taf} by restricting to abelian varieties with complex multiplication by a fixed imaginary quadratic field and a compatible polarization; the corresponding Shimura stack yields the spectrum of topological automorphic forms.

Some combination of these ideas could produce a spectrum analogous to $\mathrm{TMF}$ for K3 surfaces: a ``universal K3 cohomology theory.'' Is there a way to define modular forms for K3 surfaces, so that this conjectural spectrum constitutes ``topological modular forms'' just as $\mathrm{TMF}$ does for elliptic curves? If this spectrum can be constructed, then to what extent is it oriented? The analogous construction on the height-$1$ complex $K$-theory spectrum yields a $\mathrm{Spin}$-oriented spectrum; $\mathrm{TMF}$ is $\mathrm{String}$-oriented; and so presumably this spectrum would possess some even weaker orientation.

\subsection{An Arithmetic Interpretation of the Gamma Family}

Let $p$ be an odd prime, $\BP$ denote the Brown--Peterson spectrum at $p$, and $\pi_k^S$ denote the $k$th stable homotopy group of spheres. The Adams--Novikov spectral sequence
\begin{equation*}
E_2^{s,t}=\Ext_{\BP_*\BP}^{s,t}(\BP_*,\BP_*)\implies(\pi_{t-s}^S)_{(p)}
\end{equation*}
is an important tool for explaining periodicities in the $p$-local stable homotopy groups of spheres. Its structure increasingly exhibits deep connections to arithmetic. For instance, its $1$-line is generated by elements $\alpha_t\in E_2^{1,2(p-1)t}$, for $t\ge 1$. These can be placed in correspondence with certain Bernoulli numbers: the order of $\alpha_t$ is equal to the $p$-factor of the denominator of the quotient $B_{(p-1)t}/(p-1)t$ (in lowest terms). Similarly, when $p\ge 5$, its $2$-line contains nontrivial elements $\beta_t\in E_2^{2,*}$, for each $t\ge 1$. Behrens \cite{behrens} has placed these $\beta_t$ in correspondence with integral modular forms for the congruence subgroup $\Gamma_0(1)$ of level $(p^2-1)t$ satisfying certain conditions on their $q$-expansions.

Now, for $p\ge 7$, there is also a family of nontrivial elements $\gamma_t\in E_2^{3,*}$ in the $3$-line of this spectral sequence. No arithmetic interpretation of these elements is known; however, since the Bernoulli numbers arise in connection with the height-$1$ complex $K$-theory spectrum, and modular forms arise in connection with height-$2$ elliptic spectra and $\mathrm{TMF}$, the series $\gamma_t$ might admit a similar description in connection with height-$3$ K3 spectra. We hope to find such a description in terms of arithmetic invariants of K3 surfaces, perhaps using the conjectural ``K3 modular forms'' of \S\ref{sec:tmf}.

\bibliographystyle{merlin}
\bibliography{references}

\end{document}